\documentclass[12pt]{article}
\usepackage{graphicx}
\usepackage[latin1]{inputenc}
\usepackage{amsmath,amssymb}
\usepackage{amsmath}
\usepackage{amsfonts}
\usepackage{amssymb}
\begin{document}

\title{Pseudospheres in geometry and physics: from Beltrami to De Sitter and beyond \footnote{\small Based on a lecture presented at the meeting: {\it Un grande matematico dell '800: omaggio a Eugenio Beltrami}, Milano 
14-15 October 2004. To be published in the Proceedings of the meeting by the Istituto Lombardo di Scienze e Lettere.} }

\author{Bruno Bertotti\footnote{
Dipartimento di Fisica Nucleare e Teorica, Università degli studi di Pavia.} ,  Roberto Catenacci\footnote{ Dipartimento di Scienze e Tecnologie Avanzate, Universit\`{a} del Piemonte Orientale.}  and Claudio Dappiaggi \footnote{Dipartimento di Fisica Nucleare e Teorica \& Istituto Nazionale di Fisica Nucleare, sezione di Pavia, Università degli studi di Pavia.}}

\date{}
\maketitle

\begin{abstract}
We review the extraordinary fertility and proliferation in mathematics and physics of the concept of a surface with constant and negative Gaussian curvature. In his outstanding 1868 paper Beltrami discussed how non-Euclidean geometry is actually realized and displayed in a disk on the plane. This metric is intrinsically defined and definite; but only if indefinite metrics are introduced it is possible to fully understand the structure of pseudospheres. In a three-dimensional flat space $\mathbb{R}^{3}$ the {\em fundamental quadric} is introduced, with the same signature as the metric of $\mathbb{R}^{3}$; this leads to three kinds of surfaces with constant Gaussian curvature: the sphere, the single-sheet hyperboloid and the two-sheet hyperboloid; the last one is shown to be isomorphic to Beltrami's disk. If a hyperboloid is extended to a four-dimensional metric with signature $(+---)$ (as in spacetime), the two-sheet case describes relativistic kinematics of free particles, an example of non Euclidean geometry already recognized in 1910. The spacetime corresponding to the single-sheet case is {\em de Sitter cosmological model}, which, due to its symmetry, has an important role in cosmology. When two of the three fundamental quadrics are combined, a simple, yet deep, solution of Einstein-Maxwell equations corresponding to a uniform electromagnetic field is obtained (\cite{bb59b, robinson59}, here called BR; see also \cite{Levi-Civita}) with many applications in mathematical physics. One of them is a `no go' theorem: when one tries to frame a Riemannian four-dimensional manifold in a K\"{a}hlerian structure, it is found that, while this is generically possible with a definite signature, in spacetime only BR fulfills the requirement. The BR metric plays an important role in the exploration of new principles in fundamental physical theories; we discuss some examples related to the horizon of a black hole and the dilaton in string theory.
\end{abstract}

\section{Beltrami and non-Euclidean geometry}

Since the axioms of Euclidean geometry, in particular, the properties of parallel lines, cannot be demonstrated, the development of a self-consistent non Euclidean geometry based upon different postulates has been for a long
time an important goal of mathematicians, with great relevance for the foundations of mathematics itself \cite{bonola06, fano58, rosenfeld88}. Until the end of the XIX century much of this work -- in particular Lobatschewsky's, Gauss' and Bolyai's -- is based upon an abstract reformulation and development of different postulates and geometrical entities which do not admit intuitive representations, but acquire their meaning in reciprocal relationships. These abstract developments are closely connected with projective geometry.

Important intuitive aspects, however, had emerged in connection with the theory of surfaces which, with their geodesics, provide a realization of `straight lines' and non Euclidean geometries; in particular, Gauss in 1827
showed that the sum of the internal angles of a geodesic triangle in general is not $\pi$;
\begin{equation}
\alpha+\beta+\gamma-\pi=\varepsilon=KA.\label{egregium}
\end{equation}
is the {\em excess angle} of the triangle, proportional to its area $A$. When $K=0,\, >0,\,<0$ we have Euclidean, {\em elliptic} and {\em hyperbolic} geometry, respectively. $K$ is the local {\em Gaussian curvature} of the surface, a quantity invariant for any isometric deformation. For example, as indicated in the leftmost diagram of Fig. 3, in the spherical triangle $ABP$, of area $\alpha R^{2}$, the excess angle is just $\alpha$. In an ordinary surface, the sections with planes through the normal at a point $P$ define a one-parameter set of curves; if $R_{_{M}}$ and $R_{m}$ are, respectively, the largest and the smallest radius of curvature, taken with their signs, $K=1/R_{_{M}}R_{m}$. A point where $K<0$ is a saddle point. Of course, spherical trigonometry (when $K=1/R^2$ is constant and positive), which was well known much earlier, provides a straightforward realization of elliptic non Euclidean geometry. It is surprising that until the middle of the XIX century this was not recognized; the term non Euclidean geometry was reserved to the hyperbolic case. As hinted by Gray (p. 156 in \cite{gray90}) this was due to the fact that the geodesics on a sphere -- the great circles -- have a finite length, in contradiction with the usual concept of an infinite `straight line'. This goes back to Saccheri \cite{saccherio33}, who claimed to show that a geometry fulfilling the hypothesis of an `obtuse angle' (equivalent to $\epsilon >0$) is contradictory; but in the lengthy discussions of Propositions XI and XII he implicitly assumes that a `straight line' is infinite (see also \cite{bonola06}). We know well that, broadly speaking, the concept of `straightness' can be realized mathematically either in an affine way, as autoparallelism, or metrically: a straight segment is the shortest line between two points. None of them requires infinite extension. 

In his masterly paper \cite{beltrami68}, Beltrami investigated a detailed realization of hyperbolic geometry in an open disk in a plane, with constant and negative Gaussian curvature $K=-1/R^{2}<0$ \footnote{For an exhaustive discussion of Beltrami's work on non Euclidean geometry, including his correspondence with J. Ho\"{u}el, see \cite{boi98}.}. It is intrinsically defined by the metric
\begin{equation}
ds_{B}^{2}=R^{2}\frac{(R^{2}-v^{2})du^{2}+2uv\,du\,dv+(R^{2}-u^{2})dv^{2}
}{(R^{2}-u^{2}-v^{2})^{2}}.\label{eq:beltrami}
\end{equation}
The variables $u$ and $v$ are confined to the open disk $u^{2}+v^{2}<R^{2}$, the points corresponding to its boundary being at infinity. In this representation geodesics are just straight segments in the $(u,v)$ plane and
the model is geodesically complete (see Fig. 1) (Beltrami, in his earlier paper \cite{beltrami65}, had already solved the general problem, how to find coordinates on a surface such that geodesics are straight lines.) In this way all the peculiarities and the trigonometry of hyperbolic geometry were clearly demonstrated with the methods of differential geometry; moreover, the door was opened to its fertile study in a projective framework.

\begin{figure}[h]
\begin{center}
\includegraphics[width=1.5in]{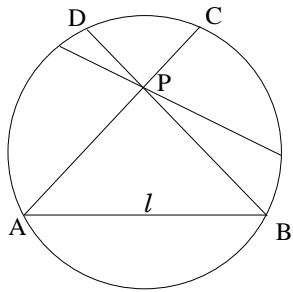}
\end{center}
\vspace*{0pt}
\caption{Beltrami's realization of a metric of constant and negative Gaussian curvature on the open disc $u^{2}+v^{2} < R^{2}$ in the plane $(u,v)$. A (straight) chord $\ell$ corresponds to a geodesic, its end points $A, B$ being at infinity. Through any point $P$ there are infinitely many geodesics parallel to $AB$, the extremal elements of this class being $CA$ and $DB$.}
\end{figure}

Is it possible to realize the geometry of Beltrami's disk in the same way as a sphere in ordinary three-dimensional space realizes the elliptic non Euclidean geometry? Earlier Minding \cite{minding39, minding40} had studied surfaces of constant Gaussian curvature and shown that, if $K<0$, they are realized by
rotating the {\em tractrix}
\begin{equation}
\xi=\frac{R}{{\cosh\chi}},\,\,\,\zeta=R(\chi-\tanh\chi) \label{eq:tractrix}%
\end{equation}
around the $\zeta$-axis (Fig. 2)\footnote{At the Department of Mathematics of the University of Pavia there are several models of this and related surfaces, including a large one built with paper by Beltrami himself (http://www-dimat.unipv.it/bisi/cuffia).}. The curvatures of the meridian and the parallel sections have opposite signs, so that each point is a saddle point. Contrary to the Beltrami disk, however, the surface so generated is neither complete, nor simply connected. In 1901 D. Hilbert \cite{hilbert01} proved that in Euclidean space there is no regular and complete surface where Lobatschewsky geometry can be realized; clearly this is not the way to obtain its correct representation. Minding's surface is invariant only under rotations around the $\zeta$-axis; for the sphere, instead, we have the full, three-dimensional rotation group. The fact that a similar invariance group holds for Beltrami's metric is vaguely hinted at by Beltrami in the introduction of \cite{beltrami68}, where he recalls that Gauss' discovery that surfaces of constant Gaussian curvature are `unconditionally applicable' on themselves. Beltrami came very near the solution when, in the first note at the end, observes that the metric obtained from (\ref{eq:beltrami}) by making $R=iR^{\prime}$ imaginary covers the whole plane $(u,v)$ and is generated by the projection of the points of a sphere of radius $R^{\prime}$ from its center onto a tangent plane. The step to the proper embedding appears now straightforward, but requires acceptance of, and familiarity with, indefinite metrics. How can a surface of constant and negative Gaussian curvature be embedded in a flat three-dimensional space, preserving the invariance under its full isometry group? The answer to this question is, of course, hidden in the theory of surfaces of constant Gaussian curvature and became accessible only in the framework of Riemannian geometry, introduced by Riemann in 1867 \cite{riemann67}; but the indefinite case has been really understood much later also in relation to Einstein's theory of General Relativity (e. g., \cite{eisenhart97}). We now review in simple terms the general problem of surfaces of constant Gaussian curvature; for an excellent, although not completely general discussion, see \cite{cartan63}. 

\begin{figure}[h]
\begin{center}
\includegraphics [scale=.6]{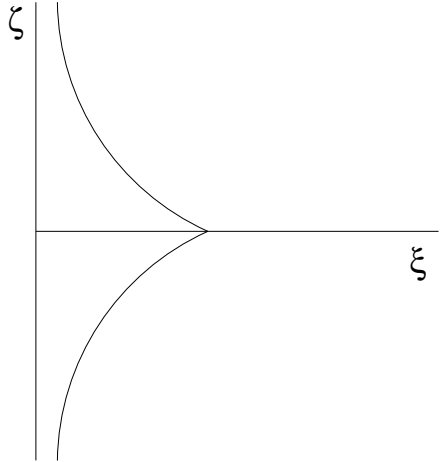}
\end{center}
\caption{The surface of negative and constant Gaussian curvature in ordinary space is obtained by rotating the tractrix (\ref{eq:tractrix}) around the $\zeta$-axis.}
\end{figure}

\section{Pseudospheres}

The requirement of a positive metric must be abandoned. In a (real) flat embedding space $(\mathbb{R}^{3},d\sigma^{2})$ with metric
\begin{equation}
d\sigma^{2}=\epsilon_{\xi}d\xi^{2}+\epsilon_{\eta}d\eta^{2}+\epsilon_{\zeta
}d\zeta^{2},\label{metric}
\end{equation}
consider the \emph{fundamental quadric} \cite{eisenhart97}
\begin{equation}
\epsilon_{\xi}\xi^{2}+\epsilon_{\eta}\eta^{2}+\epsilon_{\zeta}\zeta
^{2}=\epsilon R^{2}.\label{eq:quadric}
\end{equation}
$\epsilon_{\xi},\epsilon_{\eta},\epsilon_{\zeta}$ and $\epsilon$ take the values $+1$ or $-1$. Each case will be denoted with the symbol $(sgn\,\epsilon_{\xi} \;sgn\,\epsilon_{\eta} \;sgn\, \epsilon _{\zeta},\;sgn\,\epsilon)$ (see Fig. 3). Since the two cases $(+++,-)$ and $(---,+)$ are forbidden, with a permutation of the embedding coordinates, we can choose $\epsilon_{\zeta}=\epsilon.$

A caution about the concept of Gaussian curvature is helpful here. The usual formula $K=\frac{1}{R_{M}R_{m}}$ is obvious only for positive definite metrics; in the indefinite case, a careful definition of the radius of curvature and its sign is required. Skirting this complication, the traditional expression can be extended just by keeping in all cases the definition of half the (invariant) scalar curvature, i. e. $K=\mathcal{R}/2$ . Note also that rescaling the metric $g_{ij}$ to $\alpha g_{ij}$ with $\alpha\neq 0$ changes $\mathcal{R}$ and $K$ to $\mathcal{R}/\alpha$ and $K/\alpha$. Similarly, the {\em Theorema egregium} (\ref{egregium}) cannot be directly generalized to an indefinite or to a negative metric; the concept of angle must be carefully extended, introducing hyperbolic functions \cite{nomitzu}. For this reason, and for simplicity, from now on we shall not use the concept of angle anymore.

The manifolds (\ref{eq:quadric}) with the induced metric (\ref{metric}), are complete surfaces in $\left(\mathbb{R}^{3},d\sigma^{2}\right)$ with constant Gaussian curvature $K =\epsilon/R^{2} \neq 0$. When simply connected, they are the only ones with these properties. Their geodesics are the intersections with planes through the origin $O$. Note that the sign of $K$ is solely determined by $\epsilon$; this is also the sign of the square of the normal vector $\mathbf{n}$, and can be directly found taking a particular point, e. g. $P$ (Fig. 3). When viewed as subsets of the Euclidean $\mathbb{R}^{3}$, there are three topologically distinct kinds.

\begin{figure}[h]
\includegraphics[bb=-10 -10 550 250, scale=.7]{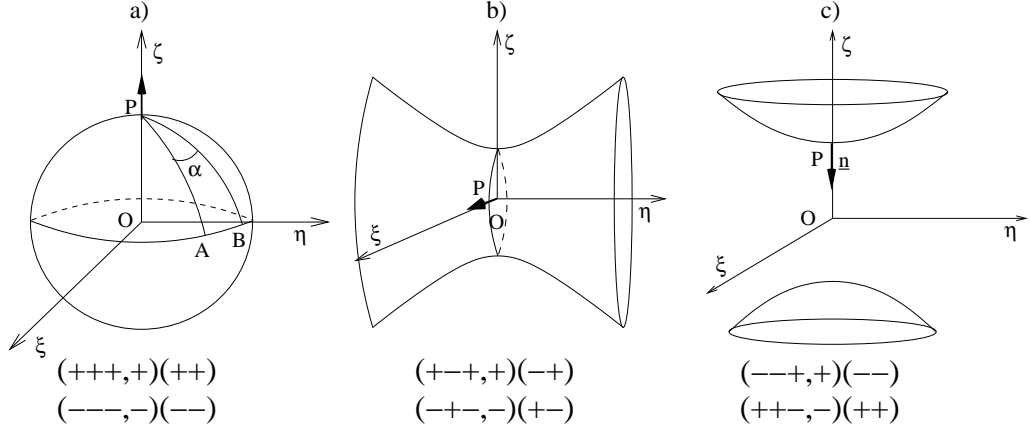}
\caption{Topological classification of the fundamental quadrics (\ref{eq:quadric}) in $\mathbb{R}^{3}$, with constant Gaussian curvature $K=\epsilon/R^{2}$. a)  The sphere $S^2$. The triangle $ABP$, of area $A=\alpha
R^{2}$, illustrates Gauss' relation (\ref{egregium}). b) The single-sheet hyperboloid. If viewed as a spacetime, it is called $dS_2$ or $AdS_2$ according to whether its axis is time-like or space-like (see Table below) c) The 
two-sheet hyperboloid $\mathbb{H}$. As explained in the text, each of these topologically different cases can be realized metrically in two ways. The direction of the vector $\mathbf{n}$ corresponds to one of them. Each realization is denoted by the signs of the four $\epsilon$'s; in the second bracket we give the intrinsic signature.}
\end{figure}

With a definite metric $(\epsilon_{\zeta}=\epsilon_{\eta}=\epsilon_{\zeta} =1)$, $\epsilon$ must also be $1$; this is the usual sphere ${S}^{2}$ (note that it is the only one that can be embedded with the induced metric in the
Euclidean $\mathbb{R}^{3}$). We have also the `negative sphere' $(---,-)$. When the metric (\ref{metric}) is indefinite, the fundamental quadric is a hyperboloid. A double-sheet hyperboloid can be realized either as $(--+,+)$ or $(++-,-)$; in both cases the induced metric $ds^{2}$ is also definite, $(--)$ in the first case and $(++)$ in the second. $(-+-,-)$ and $(+-+,+)$ are single-sheet hyperboloids; at every point there are two null lines, the intersections with the tangent plane. The induced metric is indefinite, $(-+)$ in the first case and $(+-)$ in the
second. If the coordinates $\xi$ and $\eta$ are interchanged the same is true for $(-++,+)$ and $(+--,-)$.

All these metrics appear in pairs, obviously equivalent from the point of view of the geometry of the geodesics; they differ only in the sign of scalar curvature $\mathcal{R}$. In all cases they are the locus of the points at
the same $\sigma$-distance from the origin $O$ and, therefore, they possess the three-dimensional isometry group which leaves the embedding metric invariant. Their Riemann and Ricci tensors are
\[
R_{ijkh}=K(g_{ik}g_{jh}-g_{ih}g_{jk}),\;\;R_{ij}=Kg_{ij},\;\;
\mathcal{R}=g^{ik}g^{jl}R_{ijkl}\;\;(\mathcal{R}/2=K=\epsilon/R^{2}),
\]
where Latin indices denote two intrinsic coordinates.

The two dimensional metrics can be written in a `conformally flat' form:
\[
ds^{2}=\frac{\epsilon_{\xi}du^{2}+\epsilon_{\eta}dv^{2}}{\left(
1+\frac{\epsilon}{4R^{2}}\left(  \epsilon_{\xi}u^{2}+\epsilon_{\eta}
v^{2}\right)  \right)  ^{2}},
\]
where
\[
u=\frac{2\xi}{1+\left(  1-\frac{\epsilon}{R^{2}}\left(  \epsilon_{\xi}\xi
^{2}+\epsilon_{\eta}\eta^{2}\right)  \right)  ^{\frac{1}{2}}},\qquad
v=\frac{2\eta}{1+\left(  1-\frac{\epsilon}{R^{2}}\left(  \epsilon_{\xi}\xi
^{2}+\epsilon_{\eta}\eta^{2}\right)  \right)  ^{\frac{1}{2}}},
\]
are coordinates in a `plane' with metric $dl^{2}=\epsilon_{\xi}du^{2}+\epsilon_{\eta}dv^{2}$. The pseudospheres
can also be mapped with hyperbolic variables $(-\infty<\chi <\infty,0\leq\phi<2\pi$). For the single-sheet case the intrinsic metrics are
\begin{equation}
ds^{2}=\epsilon R^{2}\left(  -d\chi^{2}+\cosh^{2}\chi d\phi^{2}\right).
\label{single-sheet}
\end{equation}
For the two-sheet case
\begin{equation}
ds^{2}=\epsilon R^{2}\left(  -d\chi^{2}-\sinh^{2}\chi d\phi^{2}\right)
.\label{double-sheet}
\end{equation}
It can easily be shown that the stereographic projection of one sheet from the origin $O$ onto the tangent plane through $P$ reproduces the metric $ds_{_{B} }^{2}$ (\ref{eq:beltrami}) or $-ds_{_{B}}^{2}$. Explicitly, the required mapping $(\chi,\phi)\rightarrow(u,v)$ is
\[
u=R\tanh\chi\cos\phi,\;\;v=R\tanh\chi\sin\phi.
\]
The points on the circumference $u^{2}+v^{2}=R^{2}$ correspond to the points at infinity on the upper null cone $\xi^{2}+\eta^{2}=\zeta^{2}$ $(\zeta>0)$. A complete representation can be obtained by formally identifying the antipodal points $(\xi,\eta,\zeta)$ and $(-\xi,-\eta,-\zeta)$.

We now understand why Minding's surface is not a good representation of Beltrami's metric (\ref{eq:beltrami}): the positive definite signature of the embedding space prevents realizing the three-parameter symmetry. Bianchi, in \cite{bianchi94}, Ch. XVI, pointed out that the symmetry group is represented in the plane $(u,v)$ by the unimodular transformations
\[
w^{\prime}=\frac{aw+b}{cw+d}\qquad(ad-bc=1),
\]
where $w=u+iv$ and $a,b,c,d$ are real numbers. It is well known that this ultimately leads to the spin representations. 

The single-sheet hyperboloid could also provide a realization of non Euclidean geometry. It is not worthwhile to dwell here in such details, but we only show that a similar stereographic projection can generate the `Beltrami
representation' of the single-sheet hyperboloid $(-+-,-)$. Denoting with $(u,v)$ the $\xi$ and $\eta$ coordinates in the projection plane $\zeta=R$, the projection $P_{0}$ of a point $P=R(\cosh\chi\cos\phi,\sinh
\chi,\cosh\chi\sin\phi)$ on the hyperboloid has $u=R\cot\phi,v=R\tanh\chi \csc\phi,$ leading to the metric:
\[
ds_{B2}^{2}=R^{2}\frac{(v^{2}-R^{2})du^{2}-2uvdudv+(u^{2}+R^{2})dv^{2}}
{(R^{2}+u^{2}-v^{2})^{2}},
\]
to be compared with the original Beltrami metric (\ref{eq:beltrami}). The straight lines in the $(u,v)$ plane are still geodesics, but Beltrami's disk is replaced by\ $v^{2}-u^{2}<R^{2}$, the points at infinity (in the appropriate
Minkowsky sense) corresponding to those on the hyperbola $v^{2}-u^{2}=R^{2}$. In this way a projective analysis of the geometry of this second Beltrami metric, as in (\ref{eq:beltrami}), is possible.

The two-dimensional surfaces generated by fundamental quadrics (Fig. 3) can be combined to construct Riemannian manifolds with an even number of dimensions endowed with interesting symmetries. For spacetime, however, the correct number of time-like (one) and space-like (three) intrinsic coordinates must be ensured, narrowing down the choice. Consider a four-dimensional manifold $\Sigma_{4}=\Sigma_{+}\otimes\Sigma_{-}$, topological product of two two-dimensional manifolds, $\Sigma_{+}$ with coordinates $(x_{0},x_{1})$, and $\Sigma_{-}$, with coordinates $(x_{2},x_{3})$; a tensor is termed \emph{decomposable} if a) its components with mixed indices vanish and b) its components relative to $\Sigma_{+}$ ($\Sigma_{-}$) depend only on the respective coordinates. A Riemannian manifold $\Sigma_4$ with this property is decomposable if its metric tensor is decomposable:
\[
g_{\mu\nu}=g_{(+)\mu\nu}+g_{(-)\mu\nu}.
\]
Its Ricci tensor has the same property:
\begin{equation}
R_{\mu\nu}=K_{+}g_{(+)\mu\nu}+K_{-}g_{(-){\mu\nu}},
\end{equation}
where $K_{+}$ and $K_{-}$ are the respective Gaussian curvatures \cite{ficken39}. Greek indices run from $0$ to $3$. As discussed in Sec. 5, they directly lead to the Bertotti-Robinson solution of Einstein-Maxwell field
equations \cite{bb59b, robinson59}.

\section{Hyperbolic geometry in special relativity}

Obviously the fundamental quadrics can be generalized to $n$-dimensional manifolds of constant Gaussian curvature $K$ embedded in a flat space $\mathbb{R}^{n+1}$.

In special relativity the energy $E$ and the momentum ${\bf p} = (p_{x},p_{y},p_{z})$ of a free particle with rest-mass $m$ fulfill\footnote{The signature $(+---)$ is assumed for spacetime, at variance with the usual choice $(-+++)$ in quantum field theory; the velocity of light is unity.}
\begin{equation}
E^{2}-p_{x}^{2}-p_{y}^{2}-p_{z}^{2}=m^{2}.\label{eq:shell}
\end{equation}
This is the extension of the two-sheet hyperboloid $(--+,+)$ to four dimensions; the two sheets correspond to particles and antiparticles, respectively. There is no wonder that the kinematics of special relativity
embodies a realization of hyperbolic non Euclidean geometry. This recognition became possible only after the concept of spacetime, as introduced by H. Minkowsky in his fundamental paper \cite{minkowsky08}, was understood and
adopted; rather than Lorentz transformations and explicit coordinates, one can just use the geometry of null cones and the orthogonality in spacetime. This applies, in particular, to the relativistic addition of velocities; its non
Euclidean significance was soon pointed out by V. Vari\'{c}ak (1865-1942) \cite{varicak12} and A. Sommerfeld (1868-1951) \cite{sommerfeld09}.

{\it A posteriori}, this recognition appears trivial. Were the momentum $\mathbf{p}$ imaginary, the mass shell (\ref{eq:shell}) would be a sphere, on which spherical trigonometry holds. A spherical triangle with a vertex $O$ at
$\mathbf{p}=0$ and two vertices $P,Q$ at the geodesic distances $\phi_{_{P}}, \phi_{_{Q}}$ fulfills
\[
\cos\phi= \cos\phi_{_{P}}\cos\phi_{_{Q}}- \sin\phi_{_{P}}\sin\phi_{_{Q}}
\cos\alpha;
\]
$\phi$ is the length of $PQ$ and $\alpha$ the angle between $OP$ and $OQ$. The corresponding formula for hyperbolic geometry is recovered when $\phi=i \chi$; then 
\[
\cos\phi= \cosh\chi= (1-v^{2})^{-1/2},
\]
and we get the usual relativistic addition formulas in terms of the \emph{celerity}, the hyperbolic angle $\chi= \tanh^{-1} v$. When the velocities are parallel ($\alpha=0$), the addition theorem just says that celerities add together.

\section{De Sitter's Universe}

De Sitter's four-dimensional manifold $dS_{4}$, generalizes the single-sheet hyperboloid $(-+-,-)$ and in a flat five-dimensional space $\mathbb{R}^{5}$ with metric 
\begin{equation}
d\sigma^{2}=-d\zeta^{2}+d\eta^{2}-d\xi^{2}-d\upsilon^{2}-d\omega^{2},
\label{eq:deSmetric}
\end{equation}
is defined by the fundamental quadric
\begin{equation}
-\zeta^{2}+\eta^{2}-\xi^{2}-\upsilon^{2}-\omega^{2}=-R^{2}.
\label{eq:deSquadric}
\end{equation}
The induced metric has the correct spacetime signature $(+---)$. This manifold is invariant under the 10-parameter isometry group consisting of the linear transformations in $\mathbb{R}^{5}$ which leave (\ref{eq:deSmetric})
invariant; when restricted to the fundamental quadric, in the limit $R\rightarrow \infty$ it corresponds
to the Lorentz group combined with the four-parameter translational group (the Poincar\'{e} group $SO(3,1)\ltimes T^{4}$)\footnote{The symbol $\ltimes$ stands for semidirect group product.} All points of $dS_4$ are equivalent.

Just like the sphere is the `perfect' surface in ordinary space (and for this reason it played such an important role in ellenistic and Tolomean cosmologies), we have here the \emph{perfect spacetime}; indeed, more
appropriately so than the sphere, because it does not just display the Universe at a single moment of time, but it encompasses its whole history, from the infinite past to the infinite future, with no time asymmetry, no
evolution and no causal chains.

For this epistemological reason it is an obvious candidate as a model of the Cosmos itself, a peculiar object of investigation indeed, since it is given to us only once: by definition, the Cosmos does not allow the characteristic pattern of theoretical conjecture and experimental repetition and confirmation. Basing a model of the Universe on the criterion of absolute simplicity, rather than on the extrapolation of ordinary physical laws, surely seems a brilliant way to escape this epistemological instability.

It is surprising that only in 1948 H. Bondi, T. Gold and F. Hoyle in their pioneer papers \cite{bondi48, hoyle48}, recognized that the de Sitter manifold could in fact provide an acceptable model of the Universe. The
\emph{Steady State Theory} was thus created (see \cite{schrodinger56} for an extensive discussion). Their great merit was to establish a connection between the abstract, single-sheet hyperboloid in $\mathbb{R}^{5}$ and the actual
expanding Universe, and to show that the ensuing very precise and stringent predictions are consistent with observations. This model depends on a single parameter, the Hubble constant $H$, with dimension $[T^{-1}]$; all the mean properties of the Universe follow. It does not have an origin, nor a history, nor a dynamics. It is a
real marvel that it has survived for more than 20 years the aggressive comparison with evermore stringent and extensive observational constraints.

Universal expansion requires the introduction in spacetime of a \emph{cosmological substratum}, a family of three-dimensional space-like simultaneity surfaces $\Sigma_{t}$: $(t =$ const); after having averaged out
local motion, in this frame galaxies are, in the average, at rest. Current scientific cosmology is based upon the \emph{Cosmological Principle}, according to which all points and all directions of these surfaces are
equivalent, leading to \emph{Robertson-Walker spacetime} 
\[
ds^{2} = dt^{2} -R^{2}(t)\frac{d{\bar x}^{2}+d{\bar y}^{2}+d{\bar z}^{2}}{[1
+k({\bar x}^{2}+{\bar y}^{2}+{\bar z}^{2})/4]^{2}}.
\]
The dimensionless coordinates $({\bar x},{\bar y},{\bar z})$ label points moving with the substratum; the spatial sections $t =$ const are spheres, or pseudospheres, when $k=1$ and $k=-1$, respectively, and infinite Euclidean
spaces when $k=0$. The  mean global properties of astronomical observations, like source counts and the redshift-luminosity relation, are determined by the dynamically determined function $R(t)$ and the number $k$. Hubble constant $H(t) = d\ln R(t)/dt$ gives the expansion scale. In the simplest case $k=0$ a radial null geodesic issuing from the origin ${\bar r} = \sqrt{{\bar x}^{2}+{\bar y}^{2}+{\bar z}^{2}} =0$ is given by
\[
{\bar r} = \int_{t_{0}}^{t_{1}}\frac{dt}{ R(t)};
\]
then the frequency shift reads
\[
\frac{\nu_{0}}{\nu_{1}}= \frac{R(t_{1})}{R(t_{0})}= \exp\left[  \int_{t_{0}%
}^{t_{1}} dt H(t)\right].
\]

In de Sitter's world there is no evolution; hence $H$ must be constant,leading to the simple metric (with $R(0)=R=1)$
\[
ds^{2} = d{\bar t}^{2} - \exp(2{\bar t})(d{\bar x}^{2}+d{\bar y}^{2}+d{\bar
z}^{2}) \qquad({\bar t} =Ht).
\]
To recover the embedded representation, it is enough to confine oneself to the section ${\bar y}={\bar z}=0$, corresponding to the reduced quadric 
\[
\eta^{2} - \xi^{2} -\zeta^{2} =-1,
\]
and the metric in $\mathbb{R}^{3}$
\[
d\sigma^{2} = d\eta^{2} -d\xi^{2} -d\zeta^{2}.
\]
This reduced de Sitter spacetime $dS_{2}$ is parametrically represented in
$\mathbb{R}^{3}$ by
\[
\eta+ \zeta= e^{\bar t}, \quad\xi= {\bar x}e^{\bar t}, \quad\eta-\zeta={\bar
x}^{2}e^{\bar t}-e^{-{\bar t}}.
\]
This, however, maps only the half $\eta+\zeta \leq 0$; identification of the antipodal points would make the map complete \cite{schrodinger56}. The points of the (one-dimensional) substratum, on which
\[
{\bar x}^{2} = \left(  \frac{\xi}{\zeta+ \eta}\right)  ^{2} =\mathrm{const},
\]
move in $\mathbb{R}^{3}$ on planes through the origin, hence along geodesics in $dS_{2}$. The space sections ${\bar t}= $ const are represented by planes inclined by $45^{\circ}$ to the $\eta$-axis and tangent to the light cone; all
their intervals are space-like, except those with $d\xi=0$, which are null. Its symmetry properties guarantees that De Sitter's world $dS_{4}$ is a
solution of
\begin{equation}
R_{\mu\nu}=\Lambda g_{\mu\nu}\quad(\Lambda=1/R^{2}).
\end{equation}
There is no energy-momentum contribution: matter and radiation have a negligible effect on the curvature, which is entirely determined by the {\em cosmological constant} $\Lambda$. Consistently with the geometrical symmetry, Bondi, Gold and Hoyle assumed that the content of matter and radiation also fulfills the \emph{Perfect Cosmological Principle}, according to which the mean observable properties of the Universe are immutable not only from place to place, but also in time. The mean matter density must be constant; to reconcile this with the exponential
expansion, new matter must be continuously created, violating the conservation law. It is easy to see that this violation is far too small to be detectable. Moreover, the mean thermodynamical and radiative state of the Universe
does not evolve at all. 

The first evidence against the Steady State Theory came in the late 50's with the counts of radio sources, which, if
extragalactic, show clear evidence of an evolution in the distant past \cite{sciama63}; but, of course, it was the discovery of the microwave background radiation in 1965 that imposed, without any lingering doubt, the
current Big Bang paradigm (For a historical review, see  \cite{bbrbsbam90}.)

If in the fundamental quadric (\ref{eq:deSquadric}) $\eta$ and $\zeta$ are exchanged, we obtain the four-dimensional anti-de Sitter's Universe, in which the time $t$ is cyclic: every physical field defined is a periodic function of $t$. With obvious generalizations, $dS_{d}$ is the de Sitter Universe in $d$ dimensions (embedded in $\mathbb{R}^{d+1}$), one of which is time-like and open; similarly, in the $AdS_{d}$ universe the time is cyclic. 

\section{Bertotti-Robinson solution and the \emph{Already Unified Theory}}

With the geometrical tools discussed in Sec. 2, it is easy to get the physical interpretation of the topological product of two fundamental quadrics in the framework of Maxwell-Einstein theory. A (covariantly) constant and
not null\footnote{The null case is degenerate and it will not be discussed here.} electromagnetic field determines the traceless energy momentum tensor $\tau_{\mu\nu}$ in terms of its energy density $\rho$ (with dimensions
$L^{-2}$); if, in addition, we have a cosmological constant $\Lambda$, the field equations read:
\begin{equation}
R_{\mu\nu}=\tau_{\mu\nu}+\Lambda g_{\mu\nu}.\label{eq:einsteinmaxwell}
\end{equation}
The two-parameter Bertotti-Robinson solution can be obtained in a purely geometrical way, without solving any differential equation \cite{bb59a, bb59b}; Robinson has investigated the case with vanishing cosmological
constant \cite{robinson59}, while earlier Kasner \cite{kasner25} found the solution without electromagnetic field. In the wake of Einstein's and Schr\"odinger's attempts to build a unified theory of gravity and electromagnetism, the
`Already Unified Theory' \cite{misner57, rainich25} was developed, based upon the principle that geometry is all;  Einstein-Maxwell equations can be formulated in a purely geometric form, which allow to extract the
electromagnetic field from its imprint on the metric. There are three questions:

\begin{enumerate}
\item Does a solution of Maxwell-Einstein (\ref{eq:einsteinmaxwell}) fulfill purely geometric conditions?
\item Are they sufficient to determine it? 
\item Can the electromagnetic field be recovered?
\end{enumerate}

Spacetime is a connected manifold endowed with a local metric $g$. In view of the developments of the following Section, we now use the forms language and adopt four independent local null complex
one-forms:
\begin{eqnarray*}
\theta^{0}  =n_{\mu}dx^{\mu},\quad \theta^{3}    =l_{\mu}dx^{\mu},\\ 
\theta^{1}  =-\bar{m}_{\mu}dx^{\mu},\quad \theta^{2}    =-m_{\mu}dx^{\mu}.
\end{eqnarray*}
The null fields $n_{\mu}$ and $l_{\mu}$ (real) and $m_{\mu}$ (complex) fulfill:
\begin{eqnarray*}
g(l,n)    =g_{\mu\nu}l^{\mu}n^{\nu}=1,\\
g(m,\bar{m})    =g_{\mu\nu}m^{\mu}\bar{m}^{\nu}=-1.
\end{eqnarray*}
In this null tetrad the metric tensor reads:
\[
g=\theta^{0}\otimes\theta^{3}+\theta^{3}\otimes\theta^{0}-\theta^{1}
\otimes\theta^{2}-\theta^{2}\otimes\theta^{1}.
\]
The dual of a (possibly complex) two-form $F=F_{ij}\theta^{i}\wedge\theta^{j}$ is defined as\footnote{In this section Latin indices from $i$ on range from $0$ to $3$ and they label the components in the null tetrad. For an arbitrary
frame, one has $\ast F_{\mu\nu}=\sqrt{-\det g}\epsilon_{\mu\nu \rho\sigma}F^{\rho\sigma}/2.$}:
\[
\ast F_{ij}=\frac{i}{2}\epsilon_{ijpq}F^{pq}.
\]
This is an anti-idempotent operator, with $\ast\ast F=-F$. A two-form $F$ is called self-dual if $\ast F=iF;$ hence
\[
\tilde{F}=\frac{1}{2}(F-i\ast F),
\]
is the self-dual part of $F.$

The complex-valued two-forms:
\begin{eqnarray*}
Z^{1} =2\sqrt{2}\theta^{0}\wedge\theta^{1},\\
Z^{2} =2\sqrt{2}\theta^{2}\wedge\theta^{3},\\
Z^{3} =2\left(  \theta^{1}\wedge\theta^{2}-\theta^{0}\wedge\theta
^{3}\right).
\end{eqnarray*}
provide a basis for the space $\Lambda_{SD}^{2}$ of the self-dual complex two-forms $\tilde{F}.$ In the space of complex two-forms two bilinear forms are defined:
\begin{equation*}
\left(  F,G\right)=\frac{1}{4}F_{ij}G^{ij},\quad 
\left\{  F,G\right\}=-\frac{1}{2}\left(  F_{ij}G_{k}^{j}+\ast F_{ij}\ast
G_{k}^{j}\right)  \theta^{i}\otimes\theta^{k}.
\end{equation*}
$\left(  F,G\right)  $ is a scalar, which defines, when restricted to $\Lambda_{SD}^{2}$ a metric and a covariant basis:
\[
Z_{1}=Z^{2},Z_{2}=Z^{1},Z_{3}=-Z^{3}.
\]
$\left\{F,G\right\}$ is a complex, traceless and symmetric tensor.

In this formalism the Cartan structure equations read:
\[
dZ^{a}=\epsilon^{abc}\sigma_{b}\wedge Z_{c}.\qquad(a,b,c=1,2,3).
\]
The three one-forms $\sigma_{b}$ are the components of the connection $\omega$ acting on a vector field $X$:
\[
\omega(X)=\nabla_{X}\theta_{i}\otimes\theta^{i}=\frac{1}{2}\left(  \sigma
_{a}(X)Z^{a}+\overline{\sigma_{a}(X)Z^{a}}\right).
\]

The source-free Maxwell equations and the electromagnetic energy-momentum tensor $\tau$ are:
\begin{equation*}
d\tilde{F}=0,\qquad \tau   =\{\tilde{F},\overline{\tilde{F}}\}.
\end{equation*}
The electromagnetic case corresponds to a real $F_{\mu\nu}$; in this case, the mathematical and physical content of Maxwell theory is fully encoded in its (complex) self-dual part. In particular, $(\tilde{F},\tilde{F})$ is the only
invariant\footnote{This is a complex invariant; its real and imaginary parts correspond to the usual real invariants $E^{2}-B^{2}$ and $\mathbf{E\cdot B}$.} one can build with $\tilde{F};$ when $\left(  \tilde{F},\tilde{F}\right)
\neq0$ the field is said non-null, and there exists a frame in which:
\[
\tilde{F}=\frac{\sqrt{2}}{2}AZ^{3},\tau=A\bar{A}\{Z^{3},\bar{Z}^{3}\},
\]
$A$ is a differentiable complex function over spacetime. The two-planes $\theta^{0}\wedge\theta^{3}$ and $\theta^{1}\wedge\theta^{2}$ are the {\em blades} of the electromagnetic field.

The first Cartan equation gives:
\[
dZ^{3}=\sigma_{1}\wedge Z^{1}-\sigma_{2}\wedge Z^{2}\equiv\psi\wedge Z^{3}.
\]
The one-form $\psi$ is geometrically important because its components are related to the well known, and widely used, Newmann-Penrose spinor coefficients \cite{rindler}. The most important result of the Already Unified Theory \cite{bb59b, Catenacci} is that the necessary and sufficient conditions for a metric to be a solution of Einstein-Maxwell equations for a non-null electromagnetic field are:

\begin{enumerate}
\item There is a frame in which the Ricci tensor is $Ricci=A\bar{A}
\{Z^{3},\bar{Z}^{3}\}$. This condition is equivalent to the `algebraic condition' of \cite{bb59b}.

\item In this frame, $d\psi=0$ (the `differential condition' of \cite{bb59b}).
\end{enumerate}

The electromagnetic self-dual two-form is then given by
\[
\tilde{F}=\frac{\sqrt{2}}{2}e^{i\alpha}AZ^{3},
\]
where $e^{i\alpha}$ is a constant duality rotation.

The metric of the topological product of two surfaces of constant curvature trivially fulfills the two conditions. For the first one, note that when $K_{+}+K_{-}=0$ the Ricci tensor is traceless; moreover
\[
Ricci=K_{+}g_{+}+K_{-}g_{-}=K_{+}g_{+}-K_{+}g_{-},
\]
and hence
\[
K_{+}\left(  \theta^{0}\otimes\theta^{3}+\theta^{3}\otimes\theta^{0}\right)
-K_{+}\left(  \theta^{1}\otimes\theta^{2}+\theta^{2}\otimes\theta^{1}\right)
=-K_{+}\{Z^{3},\bar{Z}^{3}\}.
\]
For the second condition, we obtain, since the metric has constant curvature, $\psi=0)$.

The correct signature of spacetime and, as explained below, the positive sign of the electromagnetic energy restricts the choice among the six quadrics of Fig. 3 to just two possibilities, as shown in the Table (see also Fig. 1 of 
\cite{bb59a}). BR$_1$, with a cyclic time, has been extensively used for quantum field theoretical applications. In these papers, sometimes neither the sign subtleties involved in the Gaussian curvature, nor the global aspects 
have been taken into account. 

\begin{table}[ph]

{\begin{tabular}{@{}cccc@{}} 
Symbol & Quadrics & $\Sigma_+\otimes\Sigma_-$ & Time \\
 \hline\\
BR$_1$ & $(---,-)\otimes(+-+,+)$ & $S^2\otimes AdS_2$ & cyclic with period 
$2\pi R_-$\\ 
BR$_2$ & $(-+-,-)\otimes(--+,+)$ & $dS_2\otimes\mathbb{H}_2$ &
$(-\infty,\infty)$\\ 

\end{tabular} \label{ta1}}
\caption{The two possible realizations of the BR universe.}
\end{table}
The electromagnetic field self-dual form is covariantly constant and given by:
\[
\tilde{F}=\frac{\sqrt{2}}{2}\sqrt{\rho}Z^{3},\tau=\rho\{Z^{3},\bar{Z}^{3}\},
\]
where
\[
\rho=-K_{+}=K_{-}>0.
\]
In $BR_2$ we have:
$$g_{+}   =\theta^{0}\otimes\theta^{3}+\theta^{3}\otimes\theta^{0},\;\;\;
g_{-}=-\left(  \theta^{1}\otimes\theta^{2}+\theta^{2}\otimes\theta
^{1}\right).$$
The null tetrads are:
\begin{equation}
\left\{
\begin{array}
[c]{c}
\sqrt{2}\theta^{0}=\left(  1-\frac{x^{2}}{R_{+}^{2}}\right)  ^{\frac{1}{2}
}dt+\left(  1-\frac{x^{2}}{R_{+}^{2}}\right)  ^{-\frac{1}{2}}dx,\\
\sqrt{2}\theta^{1}=i\left(  1+\frac{z^{2}}{R_{-}^{2}}\right)  ^{\frac{1}{2}
}dy+\left(  1+\frac{z^{2}}{R_{-}^{2}}\right)  ^{-\frac{1}{2}}dz,\\
\sqrt{2}\theta^{2}=-i\left(  1+\frac{z^{2}}{R_{-}^{2}}\right)  ^{\frac{1}{2}
}dy+\left(  1+\frac{z^{2}}{R_{-}^{2}}\right)  ^{-\frac{1}{2}}dz,\\
\sqrt{2}\theta^{3}=\left(  1-\frac{x^{2}}{R_{+}^{2}}\right)  ^{\frac{1}{2}
}dt-\left(  1-\frac{x^{2}}{R_{+}^{2}}\right)  ^{-\frac{1}{2}}dx.
\end{array}
\right.  \label{tetrads1}
\end{equation}
Here $x,y,z,t$ take their usual cartesian meaning when $R_{+},R_{-} \rightarrow\infty$ (see eq. (\ref{single-sheet}) and (\ref{double-sheet}) expressed in terms of the appropriate pseudospherical coordinates).

For $BR_1$ we have
$$
g_{+}=-\left(  \theta^{1}\otimes\theta^{2}+\theta^{2}\otimes\theta
^{1}\right),\;\;\;
g_{-}=\theta^{0}\otimes\theta^{3}+\theta^{3}\otimes\theta^{0},$$
with null tetrads
\begin{equation}
\left\{
\begin{array}[c]{c}
\sqrt{2}\theta^{0}=\left(  1+\frac{x^{2}}{R_{-}^{2}}\right)  ^{\frac{1}{2}
}dt+\left(  1+\frac{x^{2}}{R_{-}^{2}}\right)  ^{-\frac{1}{2}}dx,\\
\sqrt{2}\theta^{1}=i\left(  1-\frac{z^{2}}{R_{+}^{2}}\right)  ^{\frac{1}{2}
}dy+\left(  1-\frac{z^{2}}{R_{+}^{2}}\right)  ^{-\frac{1}{2}}dz,\\
\sqrt{2}\theta^{2}=-i\left(  1-\frac{z^{2}}{R_{+}^{2}}\right)  ^{\frac{1}{2}
}dy+\left(  1-\frac{z^{2}}{R_{+}^{2}}\right)  ^{-\frac{1}{2}}dz,\\
\sqrt{2}\theta^{3}=\left(  1+\frac{x^{2}}{R_{-}^{2}}\right)  ^{\frac{1}{2}
}dt-\left(  1+\frac{x^{2}}{R_{-}^{2}}\right)  ^{-\frac{1}{2}}dx,
\end{array}
\right.  \label{tetrads2}
\end{equation}
i. e.
\begin{equation}
ds^{2}=\left(  1+\frac{x^{2}}{R_{-}^{2}}\right)  dt^{2}-\left(  1+\frac{x^{2}
}{R_{-}^{2}}\right)  ^{-1}dx^{2}-\left(  1-\frac{z^{2}}{R_{+}^{2}}\right)
dy^{2}-\left(  1-\frac{z^{2}}{R_{+}^{2}}\right)  ^{-1}dz^{2}. \label{eq:BR+}
\end{equation}
In terms of pseudo-spherical coordinates $x=R_{-}\sinh\chi, t=R_{-}\phi$. In both cases, the Gaussian curvatures are $K_{+}=-R_{+}^{-2}$,$K_{-}=R_{-}^{-2}$ with $R_{+}^{2}=R_{-}^{2}$.

We can also introduce the cosmological constant $\Lambda$ by setting,
instead:
$$K_{+}=\Lambda-\rho,\;\;\;K_{-}=\Lambda+\rho,$$
so that
\[
Ricci=\tau+\Lambda g.
\]
Then we can have $R_{+}^{2}\neq R_{-}^{2}$; the positive energy condition requires $K_{-}>\Lambda$, $K_{+}<\Lambda,$ $K_{-}>K_{+}.$ An analogous discussion for the choice of the fundamental quadrics easily follows.

\section{K\"{a}hler geometry and the BR universe.}

K\"{a}hlerian geometry is a powerful tool to generate and investigate new structures connected with Riemannian manifolds. Does it also open the way to new concepts in the foundations of physics?

In the positive definite case there are many K\"{a}hler manifolds solutions of Einstein's equations, even with an electromagnetic field. They represent {\em  gravitational instantons} and are very important in quantum gravity (for example, $S^{2}\times S^{2}$)  \cite{Catenacci2}. But we shall see that for a spacetime metric there is no non-trivial K\"{a}hlerian solution of Einstein's equation; moreover, the only K\"{a}hlerian Lorentzian solution of Einstein's Maxwell equations is the Bertotti-Robinson spacetime. This shows that the methods and the concept of complex geometry can only be applied either to complex metrics (for which the signature is meaningless) or
to definite metrics (as for Euclidean gravity).

Euclidean gravity studies Einstein's equations in a Riemannian four-dimensional space with positive definite metric. In some cases the adoption of an imaginary time coordinate, obtained with a Wick rotation $t\to it$, transforms a Euclidean gravity metric in a solution of Einstein's equations with the physically correct signature. The methods of K\"{a}hlerian geometry provide to Euclidean gravity mathematical tools of exceptional effectiveness, e. g., in the definition and evaluation of partition functions and functional integrals, leading to a deep understanding of quantum gravity in a Euclidean framework. An important application of Wick's rotation is the quantum origin of the Universe, using an action principle in Euclidean space.

Complex geometry is an old and very important tool in differential geometry. For two-dimensional manifolds of definite signature, if $z=x+iy,\bar{z}=x-iy$ are complex coordinates in the plane, it is known that
the metrics of the ordinary sphere $(++)$ and the two-sheet hyperboloid $(++)$ (corresponding to each sign in the formula below) show the important property of conformal flatness:
\[
ds^{2}=\frac{dz d\bar{z}}{\left(  1\pm\frac{z\bar{z}}{4R^{2}}\right)  ^{2}};
\]
moreover they are complex K\"{a}hlerian. On the contrary, for a one-sheet hyperboloid we would need a `plane' with hyperbolic signature and two real conformal coordinates $w=x+y,\tilde{w}=x-y:$
\[
ds^{2}=\frac{dwd\tilde{w}}{\left(  1\pm\frac{w\tilde{w}}{4R^{2}}\right)^{2}
}.
\]
It is appropriate now to discuss a suitable modification of the concept of K\"{a}hler structure adapted to the hyperbolic signature. These modified complex structures turn out to be strictly related to the geometry of the electromagnetic field discussed in the previous section.

As usual, for a real hyperbolic metric $g$ on a manifold of dimension $2n$, an {\it almost Hermitian structure} $J$ is a tensor field $J$ such that:
$$J^{2}(X)=-X,\;\;\;g(J(X),J(Y))=g(X,Y).$$
The main difference with `usual' complex geometry is that $J$ must be complex, as clearly shown by the following example in two dimensions:
\[
{\rm if }\;\;g(X,Y)=X_{1}Y_{1}+X_{2}Y_{2},\;\;{\rm then }\;\;J=\left(
\begin{array}
[c]{cc}
0 & 1\\
-1 & 0
\end{array}
\right);
\]
\[
{\rm if }\;\;g(X,Y)=X_{1}Y_{1}-X_{2}Y_{2},\;\;{\rm then }\;\;J=\left(
\begin{array}
[c]{cc}
0 & i\\
i & 0
\end{array}
\right).
\]
The covariant tensor field
\[
\Omega(X,Y)=g(X,J(Y))
\]
is antisymmetric, and there exists a null frame such that $\Omega=iZ^{3}.$ If $d\Omega=0$ the structure $J$ is said {\em almost K\"{a}hlerian}. Of course, the most interesting case is when $J$ is integrable, in the sense that there exists a maximal atlas of local {\it complex} coordinates in which the {\it real}, four-dimensional metric reads:
\[
ds^{2}=g_{\alpha\tilde{\beta}}dx^{\alpha}dx^{\tilde{\beta}}\;\;{\rm where
}\;\;\alpha=0,1\;\;{\rm and }\;\;\tilde{\beta}=\tilde{0},\tilde{1}.
\]
The only allowed transformations between overlapping charts are those of the type
\begin{eqnarray*}
y^{0}  &  =y^{0}(x^{0},x^{1}),\;y^{1}=y^{1}(x^{0},x^{1}),\\
y^{\tilde{0}}  &  =y^{\tilde{0}}(x^{\tilde{0}},x^{\tilde{1}}),\;y^{\tilde{1}
}=y^{\tilde{1}}(x^{\tilde{0}},x^{\tilde{1}}).
\end{eqnarray*}
For the metric to be real, the following conditions must hold:
\[
x^{0}=\overline{x^{0}},\;x^{\tilde{0}}=\overline{x^{\tilde{0}}},\;x^{\tilde{1}
}=\overline{x^{1}}.
\]

For our purposes the most important result is that $\Omega$ is integrable if and only if $\sigma_1=\sigma_2=0$; then it follows that the BR metric is the unique, non-flat K\"{a}hlerian solution of Einstein-Maxwell
equations \cite{Catenacci3}. The electromagnetic self-dual form is proportional to the K\"{a}hler structure $\Omega.$
The coordinates that arise from the integrability are the `conformally flat coordinates', in which the BR metric shows its product structure: in the Lorentzian case the relations between hermitian and product structures are very strong. Recalling that a {\it product structure} $P$ for a real Riemannian metric $g$ is a real traceless tensor field $P$ such that:
$$P^{2}(X)=X\;\;\;g(P(X),P(Y))=g(X,Y),$$
it is easy to verify that
\[
P=J\bar{J}(X)
\]
gives a perfect equivalence between Hermitian and product structure. In the null tetrad formalism $P=\{Z^{3},\bar{Z}^{3}\}.$ This equivalence does not occur in the case of definite signature.

\section{The BR metric as a playground for high energy physics}
We have seen that the legacy of Beltrami's work on pseudospheres ranges far beyond pure mathematics. Coming to the present, we now discuss a few examples of the wide, diverse and unexpected role that the BR metric plays in the search for progress and unification in fundamental physics. As discussed by Bousso \cite{bousso02}, this search is driven by new principles -- like the Equivalence Principle for General Relativity -- which are guidelines for the construction of new theories, though these principles must first be explored and then tested. Three general areas of investigation stand out: 
\begin{itemize}
\item A black hole is not really `black', but emits thermal radiation at a rate proportional to the area of its horizon; this has prompted deep analyses of the structure of the horizon and its neighborhood. 
\item String theory, in which the fundamental entities subject to quantization are geometrical objects with a time-like and a space-like dimension, just like an ordinary string in spacetime, has emerged as a possible unified scheme of all physical interactions. Strings are embedded in an abstract, $d$-dimensional Lorentzian manifold; the  $d-4$ dimensions beyond four-dimensional spacetime are usually assumed to be a compact submanifold with a small volume, which can be probed only at exceedingly high energies. A common feature of string theory is the universal appearance of a scalar field $\phi$ -- the {\em dilaton} -- inextricably coupled with the metric.
\item The fact that the area of the horizon of a black hole has the behavior of an entropy has led to the {\em holographic principle}, according to which, at a fundamental level, the information content of a physical system is not the sum of the content of its three-dimensional parts, but resides on a two-dimensional boundary (see, e. g. \cite{susskind95} and \cite{aharoni00}). This principle may usher a revolution and a unification in our conceptions of gravity, quantum theory and particle physics at very high energies. 
\end{itemize}
The BR universes, in particular $BR_1$, based upon $AdS_2$, have played an important role in this exploration and gave rise to a very large number of papers. In the following we only mention some examples concerning the first area. In the search for realizations of holography, an important role is played by the metric $AdS_d \otimes S^d$, which reduces to $BR_1$ when $d=2$ (see \cite{Cadoni} and references therein). Since $AdS_d$ is conformally flat, it may be possible to establish a correspondence between the content of a conformally-invariant dynamics in the bulk of the manifold with its content at infinity.    

\subsection{Near-horizon limit of the Reissner-Nordstr\o m black hole}

In a global coordinate frame $(t,r^\prime,\theta,\varphi)$, with $d^2\Omega = \sin\theta\, d\theta\, d\varphi$, let us consider an extremal Reissner-Nordstr\o m spacetime (see \cite{Wald}) with a point charge $q$ and a point mass $M =|q|$:
\begin{equation}\label{RNE}
ds^2=\left(1-\frac{M}{r^\prime}\right)^2\,dt^2- \left(1-\frac{M}{r^\prime}\right)^{-2}dr^{\prime 2}-r^{\prime 2} d^2\Omega.
\end{equation}
The structure of a quantum field theory in a given spacetime has been extensively studied with this metric as a background. The main issue concerns the properties at the quantum level of a physical system near the horizon where, eventually, semiclassical analysis fails. This corresponds to the limit $r= r^\prime-M\rightarrow 0$; in this limit, and with a constant conformal rescaling, (\ref{RNE}) becomes:
\begin{equation}\label{RNBR}
ds^2=M^2\left(\frac{r^2}{M^4} dt^2- \frac{dr^2}{r^2}-d^2\Omega\right). 
\end{equation}
We now let $r$ vary over $(0,\infty)$. It can be shown that, as long as $r^2t^2 > M^4$, a coordinate transformation can be found (see \cite{Carter} and \cite{Cadoni4}) that carries it into the metric 
\begin{equation}\label{BR+}
ds^2=\left(1+\frac{x^2}{M^2}\right)dt^2- \left(1+\frac{x^2}{M^2}\right)^{-2}dx^2 -M^2d^2\Omega.
\end{equation}
The metric (\ref{RNBR}) is `half' the $BR_1$ solution of Einstein-Maxwell equations $(---,-)\otimes(+-+,+)$, product of a sphere $S^2$  and a single-sheet hyperboloid with cyclic time ($AdS_2$), with $R_-= R_+ =M$ (see (\ref{eq:BR+})); hence $r=0$ is just a coordinate singularity. One can say, {\em $BR_1$ embodies the near-horizon properties of (\ref{RNE})}. Rather than displaying the transformation, we start from the fundamental definition of the $(t,x)$ part of the $BR_1$ metric: it is the fundamental quadric $\xi^2 - \eta^2 + \zeta^2 =M^2$ embedded in $\mathbb{R}^3$ with 
metric $d\sigma^2 = d\xi^2 -d\eta^2 + d\zeta^2$. The required embedding reads:
\begin{eqnarray}
\xi(r,t)&=&\sqrt{-M^2+\left(\frac{rt}{M}\right)^2}\sinh\left[\frac{1}{2}\ln\left(\frac{t^2}{M^2}-\frac{M^2}{r^2}\right)\right],\nonumber\\
\eta(r,t)&=&\sqrt{-M^2+\left(\frac{rt}{M}\right)^2}\cosh\left[\frac{1}{2}\ln\left(\frac{t^2}{M^2}-\frac{M^2}{r^2}\right)\right],\nonumber \\
\zeta(r,t)&=&-\frac{rt}{M}. \label{eq:BR0}
\end{eqnarray}
The restriction $r^2t^2 > M^4$ is apparent. This curtailed spacetime (\ref{RNBR}) is labeled $BR^0$. 

The embedding 
\begin{eqnarray}
\xi(r,t)&=&\sqrt{M^2+r^2}\sin\left(\frac{t}{M}\right),\nonumber \\
\eta(r,t)&=&r,\nonumber \\
\zeta(r,t)&=&\sqrt{M^2+r^2}\cos\left(\frac{t}{M}\right).
\end{eqnarray}
has no restriction, and produces the full $BR_1$ metric (\ref{BR+}), now called $BR^+$. Finally, the embedding
\begin{eqnarray}
\xi(r,t)&=&\sqrt{-M^2+r^2}\sinh\left(\frac{t}{M}\right),\nonumber\\
\eta(r,t)&=&\sqrt{-M^2+r^2}\cosh\left(\frac{t}{M}\right),\nonumber \\
\zeta(r,t)&=&r
\end{eqnarray}
also reproduces $BR_1$, but with the restriction $r^2 > M^2$; this is called $BR^-$. Its metric reads 
\begin{equation}\label{RNBRE}
ds^2=\left(-1+\frac{r^2}{M^2}\right)dt^2- \left(-1+ \frac{r^2}{M^2}\right)^{-1}dr^2-M^2d^2\Omega.
\end{equation}
It should be noted that, for simplicity, in the three cases we have used the same symbols $(t,r)$ for the time and the radial coordinate, although they label different points on the quadric and do not have the same range. The entire $BR_1$ spacetime can be interpreted as the geodesic completion of $BR^0$ (\ref{RNBR}) or of $BR^-$ (\ref{RNBRE}). The three labels $^0,\; ^+$ and $^-$ refer to the sign occurring in the expressions of $g_{00}$ and $g_{11}$. Due to the presence of a Killing horizon at $r = M$, $BR^-$  is physically quite appealing.

The global structure of these spacetimes is best understood with Penrose diagrams (see \cite{Carter} for an excellent introduction). Since a two-dimensional spacetime is always conformally flat, its metric can be put in 
the form 
$$ds^2 = C(u,v)du\,dv, $$
where $u$ and $v$ are two null coordinates, chosen in such a way as to be finite at spatial infinity. It is possible, therefore, to represent an infinite spacetime in a finite sheet of paper; the $u$ and $v$ lines are conventionally drawn as straight lines at $45^\circ$. For $BR^+$ (see Fig. 4) 
\begin{equation}C(u,v) = - \left(1+ \tan^2\frac{u-v}{2}\right), \label{conformal} \end{equation}
with the mapping 
$$ x= M\tan \frac{u-v}{2}, \quad t=M\frac{u+v}{2} \quad (0<v<2\pi, \: -\pi < u <
\pi). $$
All points whose time coordinate $t$ differs by $2\pi M$ are identified.

\begin{figure}[h]
\begin{center}
\includegraphics[width=10cm]{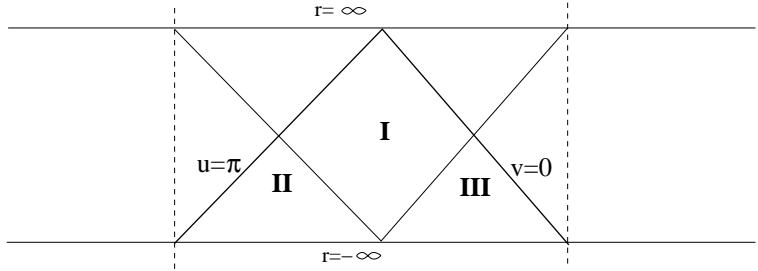}
\end{center}
\label{Penrose}
\caption{The Penrose diagram for the $AdS$ component of the $BR_1$ metric; the time $t$ runs from left to right, with period $2\pi M$, repeating indefinitely the block between the parallel lines. The $r$-coordinate is the one used for the complete $BR_1$ metric. $BR^-$ is the restriction of $BR^+$ to region {\bf I}, whereas $BR^0$ is the union of {\bf I}, {\bf II} and {\bf III} with the exclusion of the bold lines representing the locus $r=0$ in (\ref{RNBR}).}
\end{figure}
   
\subsection{Dilatonic models and the BR universe }
In the first model electromagnetism and the dilaton are coupled with gravity through the four-dimensional 
action \cite{Cadoni4}
\begin{equation}\label{DGRAV}
S=\int d^4x \sqrt{\mid g\mid} e^{-2\phi}({\cal R}-F_{\mu \nu}F^{\mu \nu});
\end{equation}
${\cal R}$ is the scalar curvature. Both the electromagnetic field $F_{\mu \nu}$ and the scalar $\phi$ act as a source for gravity. The electromagnetic Lagrangian appears with the factor $\exp(-2\phi)$, and so is the electromagnetic energy; as a consequence, the electric and magnetic binding energies of a neutral body change when it moves in a gradient of $\phi$, thereby violating the Weak Equivalence Principle. This is an example of the fact that in string theory the Equivalence Principle is generically violated. 

The spherically symmetric solution depends upon two lengths $R_+, R_-$, and reads:
\begin{eqnarray}\label{dyonic}
ds^2&=&\left(1-\frac{R_+}{r}\right)dt^2-\left(1-\frac{R_+}{r}\right)^{-1}\left(1-\frac{R_-}{r}\right)^{-1}dr^2-r^2d^2\Omega,\nonumber\\
F_{\mu\nu}&=&\frac{2q}{\sqrt{3}r^2}\epsilon_{\mu\nu\rho \sigma}u^\rho v^\sigma,\nonumber\\
\phi-\phi_0&=&-\frac{1}{4}\ln\left(1-\frac{R_+}{r}\right), 
\end{eqnarray}
where $r>R_+$, $\phi_0$ is an integration constant and 
$$2M=R_++\frac{3}{2}R_-,\;\;\;q^2=R_+R_-.$$
$q$ is the charge; $u^\mu, v^\nu$ are two unitary and orthogonal space-like vectors tangent to the surface $r=$ const. The solution has two horizons at $r= R_+$ and $r=R_-$. The extremal limit corresponds to 
$R_+-R_- =O(r-R_+) \ll 1$ and, with the variable $\eta={\rm arcsinh}\sqrt{\frac{r-R_+}{R_+-R_-}}$, 
yields the $BR^-$ component of the $BR_1$ metric in the form
\begin{equation}
ds^2=q^2(4\sinh^2\eta \,dt^2-4d\eta^2-d^2\Omega).
\end{equation}
An important geometrical difference between the extremal Reissner-Nordstr\o- m solution (\ref{RNBR}) and the metric (\ref{dyonic}) should be emphasized: in the former, the radii of curvature of the two components are equal in modulus, whereas in the latter they are arbitrary. 

The second example is the Jackiw-Teitelboim model (JT) of two-dimensio- nal gravity, with the `cosmological constant' $\Lambda$. Two-dimensional dilatonic models have important and widely ranging applications, from  toy models of quantum gravity to string theory (see \cite{Grumiller}). In the JT model gravity is coupled to the dilaton scalar $\phi$ with the action  
\begin{equation}\label{JT}
S=\frac{1}{2\pi}\int\sqrt{\mid g\mid}d^2 x e^{-2\phi}(R+ 2\Lambda).
\end{equation}
Its Euler-Lagrange equations admit the solution:
\begin{equation}
ds^2=(\Lambda r^2-a^2)dt^2-(\Lambda r^2-a^2)^{-1}dr^2, \quad \phi-\phi_0=-\frac{1}{2}\ln\left(\sqrt{\Lambda} r\right), \label{met2}
\end{equation}
where $a^2$ is a dimensionless integration constant. 
In (\ref{met2}) the metric coincides with the $(t,r)$ part of (\ref{RNBRE}): a BR universe arising in a Einstein-Maxwell-dilaton system has the same structure as two-dimensional metric produced by a dilaton. This is a key feature from the holographic point of view.

\section*{Acknowledgments}

\noindent We thank Prof. R. Tazzioli and Prof. P. L. Ferrari for suggestions about the historical part.

\end{document}